\documentclass{amsart}
\usepackage[latin1]{inputenc}
\usepackage[T1]{fontenc}
\usepackage{xtdeqnra}
\usepackage{fancyheadings}
\usepackage{amsmath}
\usepackage{amssymb}
\usepackage{amsfonts}
\usepackage{graphicx}
\usepackage{psfrag}
\usepackage{multicol}
\usepackage{setspace}
\usepackage{epsf}
\usepackage{here}

\newtheorem{prop}{Proposition}[section]
\newtheorem{thm}{Theorem}[section]
\newtheorem{defi}{Definition}[section]

\newtheorem{rmk}{Remark}[section]

\newtheorem{Hypothesis}{Hypothesis}[section]

\def\beqna{\begin{eqnarray}}
\def\eeqna{\end{eqnarray}}
\newcommand{\HH}{\mathcal{H}}

\newcommand{\R}{\mathbf{R}}

\newcommand{\Sd}{{\mathbf{S}_\delta}}
\newcommand{\Sdh}{{\mathbf{S}_{\delta^h}}}
\newcommand{\Th}{{\mathcal{T}^h}}
\DeclareMathOperator{\I}{I}
\DeclareMathOperator{\II}{II}
\DeclareMathOperator{\supp}{supp}
\DeclareMathOperator{\argmin}{argmin}
\DeclareMathOperator{\vect}{span}

\title[Numerical study of a new global minimizer for the Mumford-Shah Functional]{Numerical study of a new global minimizer for the Mumford-Shah functional in $\R^3$}
\author{Beno\^\i t Merlet}

\begin{document}
\maketitle
\begin{abstract}
In~\cite{David2}, G. David suggested a new type of global minimizer for the Mumford-Shah functional in $\R^3$, for which the singular sets belong to a three parameters family of sets ($0<\delta_1,\delta_2,\delta_3<\pi$).  
We first derive necessary conditions satisfied by global minimizers of this family. Then we are led to study the first eigenvectors of the Laplace-Beltrami operator with Neumann boundary conditions on subdomains of $\mathbf{S}^2$ with three reentrant corners. The necessary conditions are constraints on the eigenvalue and on the ratios between the singular coefficients of the associated eigenvector. 
We use numerical methods (Singular Functions Method and Moussaoui's extraction formula) to compute the eigenvalues and the singular coefficients. We conclude that there is no $(\delta_1,\delta_2,\delta_3)$ for which the necessary conditions are satisfied  and this shows that the hypothesis was wrong.
\end{abstract}

\subsection*{Keywords: } Mumford-Shah functional, Numerical analysis, Boundary value problems for second-order, elliptic equations in domains with corners.

\subsection*{AMS classification: } 35J25, 49R50, 65N38.

\section{Introduction}
\label{s1}
The Mumford-Shah functional was introduced in \cite{MuS} as a tool for image segmentation. Let $\Omega$ be a bounded open subset of $\R^n$ (the screen) and $g$ be a bounded measurable function defined on $\Omega$ (representing the image). The functional concerns pairs $(u,K)$ where $K$ is a closed subset of $\Omega$ and $u$ is a function belonging to the Sobolev space $H^1(\Omega \setminus K)$. It is defined by
\begin{eqnarray*}
J(u,K):=\HH^{n-1}(K) + \int_{\Omega \setminus K} |\nabla u|^2 +\int_{\Omega \setminus K} |u-g|^2,
\end{eqnarray*}
where $\HH^{n-1}(K)$ is the Hausdorff measure of co-dimension 1 of $K$. Let $(u,K)$ be a minimizing pair of $J$, (wich always exists~\cite{Am,DeCL}). The third term of the functional forces $u$ to be close to $g$ while, due to the second term, $u$ has slow variation on $\Omega \setminus K$. Since no regularity is assumed for $u$ across the singular set $K$, we may hope that for such a minimizer $K$ is the hyper-surface across which $g$ has great variations, i.e : the hyper-surfaces delimiting the contours of the image. 

The main difficulty arising in the theoretical study of the minimizers is the regularity of the singular set. First, let us notice that we may remove from $K$ a set of $\HH^{n-1}$ measure 0 which is not useful. Indeed, if $(u,K)$ is a minimizer, there exists a smallest closed set $K_1\subset K$ such that $u \in H^1(\Omega \setminus K_1)$. The pair $(u,K_1)$ is called a reduced minimizer of the functional. In dimension $n=2$, Mumford and Shah conjectured that if $(u,K)$ is a reduced minimizer for $J$ and $\Omega$ is bounded and smooth, $K$ is a finite union of $C^1$ arcs of curves, that may only meet by sets of three, at their ends, and with angles of $2\pi/3$.

This conjecture still resists but there exist partial results. In particular, A. Bonnet \cite{Bonnet} showed that \textit{ in the case $n=2$, every isolated connected component of $K$ is a finite union of $C^1$ curves.}

The crucial point introduced by A. Bonnet was a blow-up process, which leads to the notion of global minimizer of the Mumford-Shah functional. One way to prove the Mumford-Shah conjecture would be to get a complete description of all the global minimizers and then, if the global minimizers turned out to be simple, go back to the minimizers of the functional in a domain. The second step would be realized by proving that if a minimizer is closed to a global minimizer (which is true via blow-up) then its singular set is smooth. 

Here we are concerned with the case $n=3$. In this context, one may conjecture that the singular set of minimizer is a finite union of $C^1$ surfaces intersecting each other on a finite number of $C^1$ curves.

Let us first describe the blow-up technique. Let $(u,K)$ be a reduced minimizer of $J$, let $x\in \Omega$ and let $t>0$, we set
\begin{eqnarray*}
\Omega_{x,t}&:=& t^{-1} (\Omega-x),\\
g_{x,t}(y)&:=& t^{-1/2} g(x+ty),\qquad \forall \, y \in \Omega_{x,t},\\
K_{x,t}&:=&t^{-1} (K-x),\\
u_{x,t}(y)&:=& t^{-1/2} u(x+ty),\qquad \forall \, y \in \Omega_{x,t}.
\end{eqnarray*} 
Then $(u_{x,t},K_{x,t})$ is a minimizer of the modified functional $J_{x,t}$ in $\Omega_{x,t}$ where
\begin{eqnarray*}
J_{x,t}(v,G)&:=& \HH^{n-1}(G) + \int_{\Omega_{x,t}\setminus G} |\nabla v|^2 +t^2\int_{\Omega_{x,t}\setminus G}|v-g_{x,t}|^2.
\end{eqnarray*}
Now, let us take a sequence $(t_k)_k\downarrow 0$ and set $(u_k,K_k):=(u_{x,t_k},K_{x,t_k})$ to simplify the notations. Such a sequence is called a \textit{blow-up sequence} of $(u,K)$ at $x$. It turns out that up to extraction,  the sequence of sets $K_k$ converges to a closed subset $K_\infty$ of $\R^3$. On the other hand, since the factor $t^{-1/2}$ tends to infinity when $t$ tends to $0$, the sequence $u_k$ may not converge to a function having finite values. To overcome this difficulty, we have to subtract from $u_k$ a function which is constant on every connected component of $\R^3\setminus K_\infty$. More precisely, we have   
\begin{thm} 
\label{thm2}
There exists a closed subset $K_\infty \subset \R^n$, a function $u_\infty \in L^1_{loc}(\R^n)$ and for each connected component $V$ of $\R^n\setminus K_\infty$, constants $(\beta_{V,k})_k$ such that  up to a subsequence,
\begin{eqnarray*}
K_k&\longrightarrow K_\infty\qquad\mbox{locally for the Hausdorff distance},
u_k-\beta_{V,k} &\longrightarrow u_\infty\qquad\mbox{in } L^1_{loc}(V).
\end{eqnarray*}
Moreover, the limit pair $(u_\infty,K_\infty)$ is a reduced global minimizer of the Mumford-Shah functional in $\R^n$ (see Definition~\ref{defi1.1}).
\end{thm} 
\begin{defi}
\label{defi1.1}
Let $K$ be a closed subset of $\R^n$ and let $u \in L^1_{loc}(\R^n)$. The pair $(u,K)$ is a global minimizer of the Mumford-Shah functional in $\R^n$ if the following properties hold.
\begin{itemize}
\item For every open ball $B$ in $\R^n$, $\HH^{n-1}(K\cap B) <\infty$ and $\int_{B\setminus K}|\nabla u|^2<\infty$.
\item For every open ball $B$ in $\R^n$, for every pair $(v,L)$ which satisfies the property above and such that  
\begin{eqnarray*}\begin{array}{ll}
\mbox{a) }L\setminus B =  K\setminus B,\qquad&\mbox{b) } \quad v_{|\R^n \setminus B}=u_{|\R^n \setminus B},
\end{array}
\end{eqnarray*}
c)  if $x,y \in \R^n \setminus (B\cup K)$ belong to a same connected component of $\R^n\setminus L$, then they are also in a same connected component of $\R^n \setminus K$, \\
then
\begin{eqnarray}
\label{inegalite}
\HH^{n-1}(K \cap B)+\int_{B \setminus K} |\nabla u|^2 &\leq & \HH^{n-1}(L \cap B)+\int_{B \setminus L} |\nabla v|^2.
\end{eqnarray}
\end{itemize}
\end{defi}

From now on, we fix $n=3$. Let us list types of reduced global minimizers $(u,K)$ that are already known. For the first four types, the function $u$ is constant on each connected component of $\R^3\setminus K$.  
\begin{itemize}
\item[(i)] $K = \emptyset$. 
\item[(ii)] $K$ is a plane.
\item[(iii)]  $K$ is the union of three half planes sharing the same edge and making angles $2\pi/3$ with each other.
\item[(iv)] $K$ is the half cone spanned by the edges of a regular tetrahedron from its center. In this case $\R^3\setminus K$ has four connected components, each one being delimited by  three infinite triangular faces.      
\item[(v)] (Cracktips) $K$ is a half plane. Choosing coordinates such that $K=\{(x,0,z),\, x\geq 0,\, z \in \R \}$, the function $u$ is defined by
\begin{eqnarray*}
u(r \cos \theta,r\sin \theta,z)&=&\varepsilon \sqrt{\cfrac{2 r}{\pi}} \cos \cfrac{\theta}{2} +C ,\qquad \forall \, r>0,\; 0<\theta <2\pi,
\end{eqnarray*}
where $C$ is a constant and $\varepsilon = \pm 1$.
\end{itemize}
If this list was complete, from Theorem \ref{thm2}, every blow-up limit of a reduced minimizer should be one of the listed global minimizers. Let us now describe an example of~\cite{David2} for which this situation seems to be wrong. The whole argument is heuristic and is far from a proof. Let $R>0$ and $C>0$, the domain is the cylinder $\Omega=\{(x,y,z),\, x^2+y^2<R,\, -R <z<R$ (see Figure~1) and $g(x,y,z):=g_0(x,y)\varphi(z)$ where
\begin{eqnarray*}
g_0(r \cos \theta, r \sin \theta) &:=&\left\{ \begin{array}{ll}
C&\quad \mbox{for }0<\theta<2\pi/3,\\
0&\quad \mbox{for }2\pi/3<\theta<4\pi/3,\\
-C&\quad \mbox{for }4\pi/3<\theta<2\pi,
\end{array}
\right.
\end{eqnarray*}
and  $\varphi$ is a smooth cut-off function satisfying $0 \leq \varphi \leq 1$ and 
\begin{eqnarray*}
\varphi(z)=\left\{ \begin{array}{ll}
1&\quad \mbox{for }z\geq 1,\\
0&\quad \mbox{for }z\leq -1.
\end{array}
\right.
\end{eqnarray*}
 
Let us now consider a minimizer $(u,K)$ of the functional $J$ associated to $\Omega$ and $g$. Since $g \equiv 0$ for small $z$, one may think that, for $z_0$ close to $-R$, $K \cap \{z=z_0\} = \emptyset$. On the other hand, for $z_0$ close to $R$, since $g=g_0$, we may suppose that $K \cap \{z=z_0\}$ is close to the union of the three segments across which $g$ jumps, i.e:  $\{(r\cos \theta,r\sin \theta,z_0)\;:\;0\leq r < R,\, \theta = 2k\pi/3,\, k=0,1,2\}$. We may then expect that $K$ is the union of three regular surfaces meeting on a curve $\{ \gamma(t)\;:\; 0\leq t \leq 1\}$ satisfying $\gamma(0)=(0,0,R)$ and $\gamma(1)=(x_0,y_0,z_0)$  with $z_0 >-R$. See Figure~2 below.

\begin{figure}[h]
\begin{center}
\begin{minipage}[c]{\linewidth} 
        \begin{minipage}[t]{0.50\linewidth}
\vspace{-.75cm}
        \begin{center}
        \label{figure1}
        \hspace{4cm}
        \psfrag{m}[][]{}
        \psfrag{n}[][]{}
        \psfrag{o}[][]{$g \equiv 0$}
        \psfrag{x}[][]{$x$}
        \psfrag{y}[][]{$y$}
        \psfrag{z}[][]{$z$}
        \psfrag{p}[][]{$g \equiv C$}
        \psfrag{q}[][]{$g \equiv -C$}
        \rotatebox{0}{
        \includegraphics[scale=0.3]{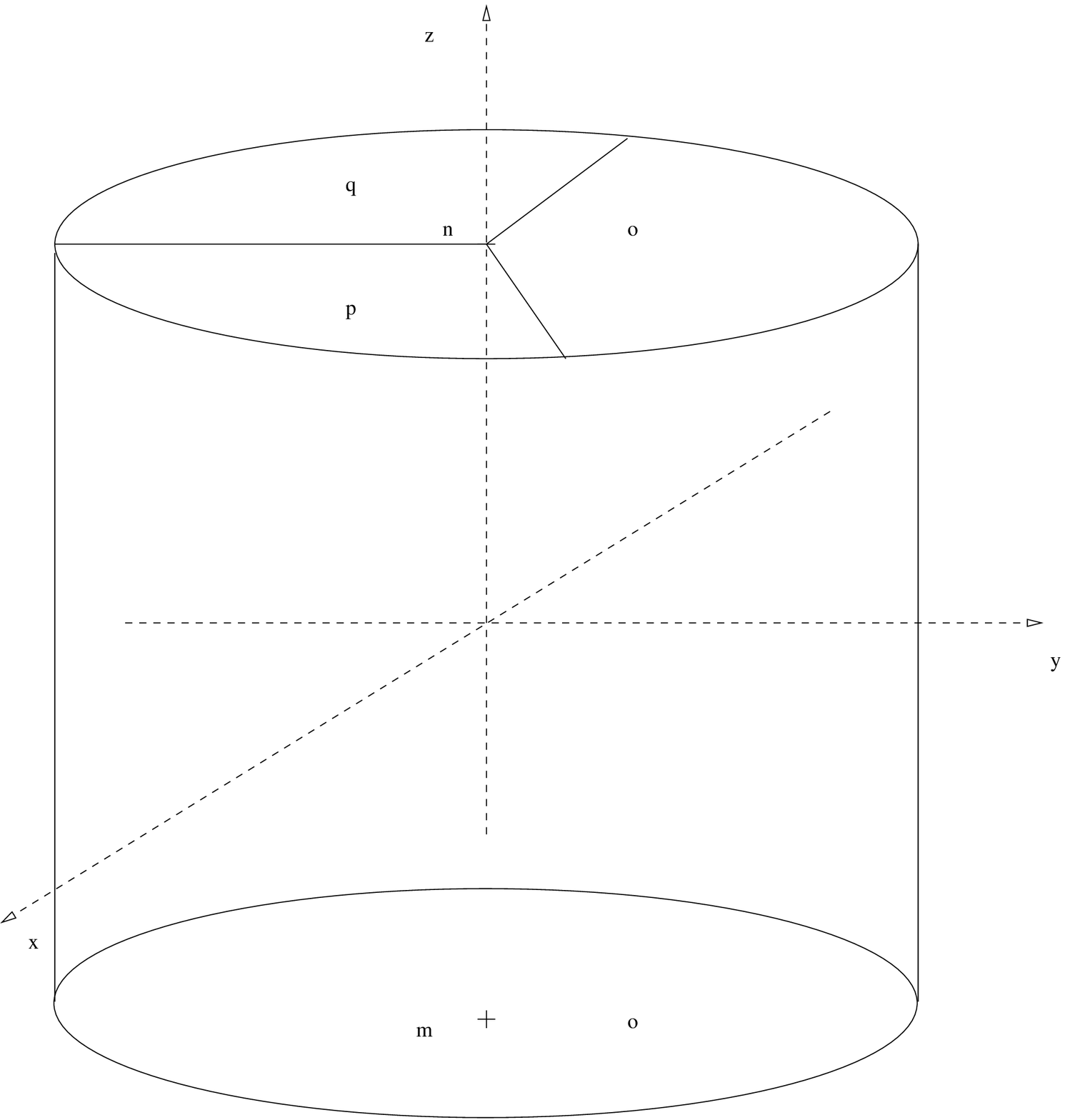}}
        \caption{Cylindrical domain $\Omega$}
        \end{center}
        \end{minipage}
        \begin{minipage}[t]{0.50\linewidth}
\vspace{0cm}
        \begin{center}
        \label{figure2}
        \rotatebox{0}{
        \includegraphics[scale=0.3]{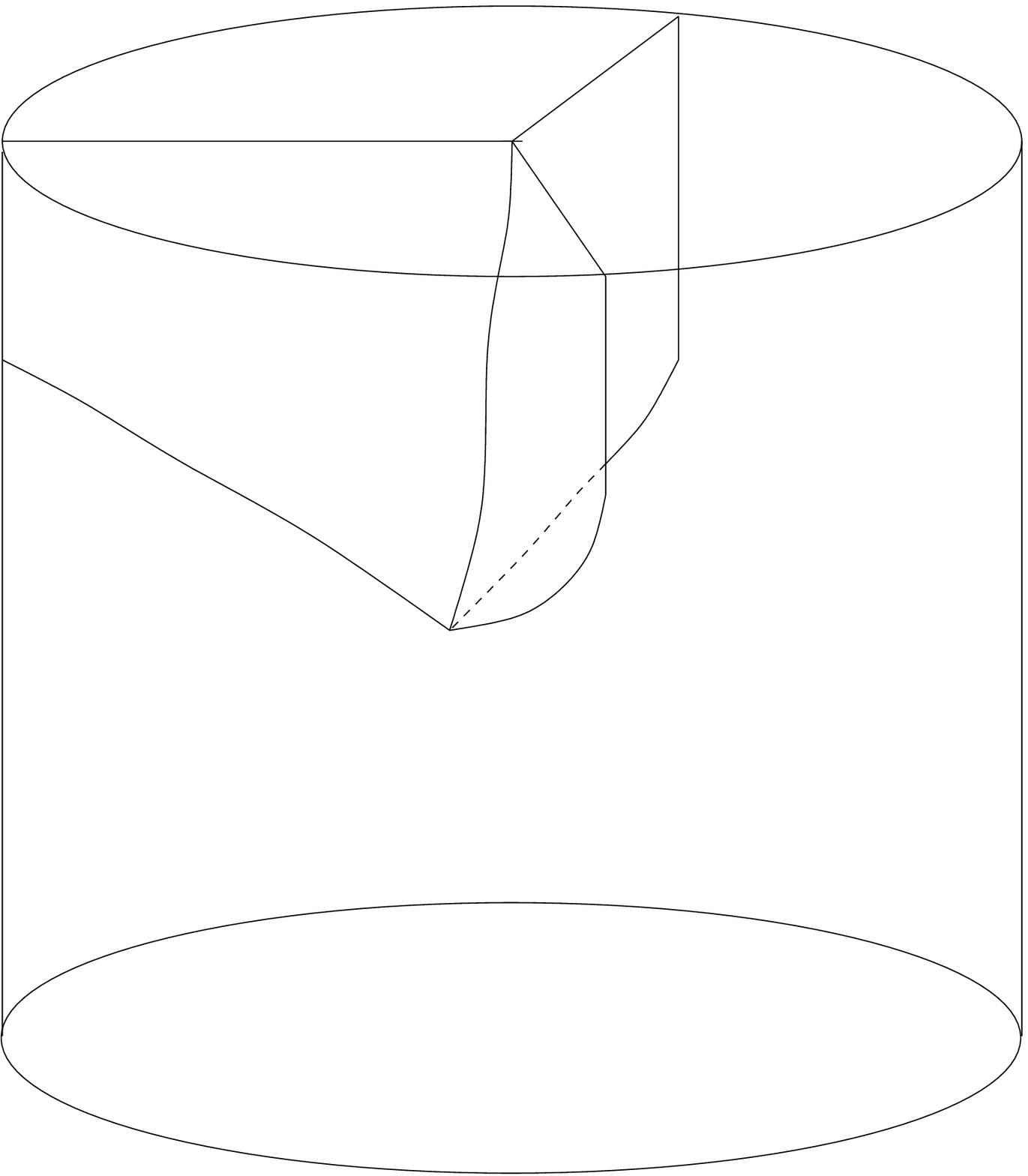}}
        \caption{Expected shape for the singular set $K$}
        \end{center}
        \end{minipage}
\end{minipage}
\end{center}
\end{figure}

For $0<t<1$, a blow-up around the point $\gamma(t)$ would lead to a global minimizer whose singular set is the union of three half planes sharing the same edge. These global minimizers should be of type (iii) and it is not necessary to introduce a new type of global minimizers if we suppose that the angles between the half plane are $2\pi/3$. The situation is different when we consider a limit blow-up at $\gamma(1)$. At this point we expect a global minimizer $(u^\star,K^\star)$ whose singular set is the union of three plane sectors with a common edge and that make angles $2\pi/3$. More precisely, intersecting $K^\star$ with the unitary sphere $\mathbf{S}^2$, we obtain a set of three arcs of big circles $M_i \subset \mathbf{S}^2$, $i=1,2,3$. These arcs are vertical, start at the north pole where they make three angles of $2\pi/3$. Denoting their lengths $\delta_i$, $i=1,2,3$, we obtain:  
\begin{eqnarray}
\label{Kstar}
K^\star = K_{\delta_1,\delta_2,\delta_3} := \R_+ \times\bigcup_{i=1}^3 M_i.
\end{eqnarray}
No global minimizer with this kind of singular set is known. And we may thing that the previous list of global minimizer was noy complete. The hypothesis of G. David is the following (\cite{David2}, sections 76 and 80):  
\begin{Hypothesis}
\label{hypothese}
There exists a new type of reduced global minimizers $(u^\star,K^\star)$ where by translation and rotation invariance there exists $0<\delta_1,\delta_2,\delta_3<\pi$ such that $K^\star= K_{\delta_1,\delta_2,\delta_3}$ and where the function $u^\star$ is homogeneous of degree $1/2$, i.e: 
\begin{eqnarray*}
u^\star(x) &=& |x|^{1/2} \Sigma^\star (x/|x|).
\end{eqnarray*} 
Moreover, this new class of global minimizers may be generated by one of them, using translation, rotation, multiplication by $-1$ and addition of a constant.

\noindent
 With this new type of global minimizers, the list of reduced global minimizers is closed.
\end{Hypothesis}
\begin{rmk}
\label{rmk1.1}
For $n=3$, the homogeneity $1/2$ is the natural homogenity for a global minimizer. In fact, if we suppose that $u^\star$ is homogeneous of degree $\alpha$ then the equilibrium between the surface term and the Dirichlet energy term in \eqref{inegalite} leads to $\alpha=1/2$ or $u$ locally constant. In our case, the homogeneity $0$ is impossible (we would remove any bounded piece of $K$ and contradict~\eqref{inegalite}). To get $\alpha=1/2$, consider $L=K\setminus B(0,r) \cup \partial B(0,r)$ and $v=0$ in $B(0,r)$, then let $R$ go to $0$ in \eqref{inegalite} to obtain $\alpha\geq 1/2$ and let $R$ go to $+\infty$ to obtain the second inequality. 
\end{rmk}

The paper is organized as follows. In section~\ref{s1.1}, we set the notations. In section~\ref{s2}, we find some necessary conditions satisfied by $(\delta_1,\delta_2,\delta_3)$ and $\Sigma$ if the hypothesis were true. Section~\ref{s3} is devoted to the description of the numerical methods we have used to check these conditions. The numerical results are presented in section~\ref{s4}. 

\subsection*{Acknowledgment}
 The author thanks Guy David for having proposed this work and for helpful support. For every fact concerning the Mumford-Shah functional, we refer to his book~\cite{David2}. We are also indebted to Patrick Ciarlet and Monique Dauge for helpful informations on the Singular Functions Method and other related methods.

\section{Notations}
\label{s1.1}
Let $\delta=(\delta_1,\delta_2,\delta_3)$ be in $(0,\pi)^3$, and let $M_1,M_2,M_3 \subset \mathbf{S}^2$ be three arcs of great circles starting from the north pole with relative angles $2\pi/3$ and with respective lengths $\delta_1,\delta_2,\delta_3$. Without loss of generality, we will assume that 
\begin{eqnarray*}
M_1 \subset \mathcal{C}_1 := \mathbf{S}^2\cap \left\{(x,y,z)\;:\; y\leq 0,\; x= 0\right\},\\
M_2 \subset \mathcal{C}_2 := \mathbf{S}^2\cap  \left\{(x,y,z)\;:\; y\geq 0, x=\sqrt{3}y \right\} ,\\
M_3 \subset \mathcal{C}_2 := \mathbf{S}^2\cap  \left\{(x,y,z)\;:\; y\geq 0, x=-\sqrt{3}y \right\}.
\end{eqnarray*}
Let $P$ be the plane $\{(x,y,z)\in \R^3\;:\; x=0\}$, in particular $\mathcal{C}_1 \subset P$. The open subset of $\mathbf{S}^2$ : $\mathbf{S}^2\setminus \cup_{i=1}^3\mathcal{C}_i$ has three connected subdomains $\Omega_{-1},\Omega_0,\Omega_{1}$. $\Omega_0$ is symmetric with respect to $P$ (i.e: $(0,1,0) \in \Omega_0$) and $\Omega_i=R^i(\Omega_0)$, $i=-1,1$, where $R$ denotes the rotation of angle $2\pi/3$ around the $z$-axis (with the usual orientation). 

We will denote by $\mathbf{S}_\delta$ the domain
\begin{eqnarray*}
\mathbf{S}_\delta &:=& \mathbf{S}^2 \setminus \bigcup_{i=1}^3 M_i.
\end{eqnarray*}
In the sequel, $L^2(\mathbf{S}_\delta)$, $H^1(\mathbf{S}_\delta)$ and $H^2(\mathbf{S}_\delta)$ will denote the standard Sobolev spaces on $\mathbf{S}_\delta$ and $\Delta$ the Laplace-Beltrami operator on $\mathbf{S}^2$. We also define the closed subspaces : 
\begin{eqnarray*}
L^2_0(\mathbf{S}_\delta)&:=&\left\{\Sigma \in L^2(\mathbf{S}_\delta)\;:\; \int_{\mathbf{S}_\delta} \Sigma =0\right\},\\
V^1(\mathbf{S}_\delta) &:=&L^2_0(\mathbf{S}_\delta)\cap H^1(\mathbf{S}_\delta).
\end{eqnarray*}
We recall the following classical result:
\begin{thm}
Let $\delta \in (0,\pi)^3$. For any $f$ in $L^2_0(\Sd)$, there exists a unique $\Sigma$ in $V^1(\Sd)$ such that 
\begin{eqnarray}
\label{diric}
\left\{
\begin{array}{rcl}
-\Delta \Sigma &=& f, \qquad \mbox{in }\Sd,\\
\partial_n \Sigma &=& 0,\qquad \mbox{on } \partial \Sd.
\end{array}
\right. 
\end{eqnarray}
Equivalently, $\Sigma$ is the unique solution in $V^1(\Sd)$ of the variational problem: $\int \nabla \Sigma \cdot \nabla S = \int f S$ for every $S$ in $V^1(\Sd)$. It is also the unique minimizer in $V^1(\Sd)$ of the functional $F(S):=1/2 \int |\nabla S|^2- \int f S$. We will note $\Sigma:=-\Delta_{N,\delta}^{-1} f$ this solution.
\end{thm}
 The operator $\Delta_N^{-1}$ is a compact symmetric operator on $L^2_0(\Sd)$. We will use spectral properties of such operators. In particular, $L^2_0(\Sd)$ has an orthonormal basis of eigenvectors of $\Delta_{N,\delta}^{-1}$. We will note $\mu_1(\delta)\geq\mu_2(\delta)\geq \cdots >0$ the eigenvalues of $-\Delta_{N,\delta}^{-1}$ counting multiplicities and for $k\geq 1$, we set $\lambda_k(\delta):=1/\mu_k(\delta)$.

\noindent
Alternatively, these eigenvalues may be defined by:
\begin{eqnarray*}
\lambda_{k}(\delta)&:=& \min_{V_k} \max \left\{\int_\Sd |\nabla S^2| \;:\; S\in V_k,\, \int_\Sd S^2 =1\right\},
\end{eqnarray*}
where the minimum is taken over all $k$-dimensional subspace $V_k$ of $L^2_{0}(\Sd)$.

\noindent
When $\delta=(0,0,0)$, it is well known that $\lambda_k(0,0,0)=2$ for $k=1,2,3$ with associated eigenvectors $(x,y,z)\mapsto x,y$ or $z$. In particular:
\begin{eqnarray}
\label{starmorgan}
\int_{\mathbf{S}^2}|\nabla S|^2 &\geq & 2 \int_{\mathbf{S}^2}S^2\qquad \forall \, S \in  V^1({\mathbf{S}^2}).
\end{eqnarray}
This property will be used at the end of section~\ref{s2}.

We will need some well known facts about the splitting in regular and singular parts of solutions to the Poisson problem with Neumann boundary conditions in a domain with corners. For this theory, we refer to~\cite{Dauge},~\cite{Grisvard} or~\cite{Kondratev}.  Let us denote by $\xi_i$ the end of $M_i$ for $i=1,2,3$. The domain $\mathbf{S}_\delta$ possesses 3 re-entrant corners of angles $2\pi$ at $\xi_1,\xi_2,\xi_3$ . 
\begin{rmk}
In~\cite{Dauge},~\cite{Grisvard}, only flat domains are considered. In order to prove Theorem~\ref{thm4}, we may use local smooth maps to transform the $-\Delta$ operator on $\mathbf{S}_\delta$ in an elliptic operator with smooth coefficients on a planar domain with a cut. In fact, it seems more natural to prove Theorem~\ref{thm4} directly. The main ingredients, Green formula, trace theorems, density results and use of polar coordinates do not change when one replace the planar domain with cuts by $\mathbf{S}_\delta$. 
\end{rmk}

We begin in introducing a set of singular functions.

\begin{defi}
\label{defi2.1}
For $x\in \mathbf{S}_\delta$, let $r_i(x)$ denote the geodesic distance on $\mathbf{S}^2$ between the points $\xi_i$ and $x$. Using the usual orientation on $\mathbf{S}^2$, for $x \in \mathbf{S}_\delta$ in the neighborhood of $\xi$, $\theta_i(x)$ denotes the angle at $\xi_i$ between $M_i$ and the smallest geodesic segments $[\xi_i,x]$ (see Figure~3 below). We use $(r_i(x),\theta_i(x))$ as polar coordinates near $\xi_i$ to define
\begin{eqnarray*}
s_i(x)&:=& 2 \tan \left( \cfrac{\sqrt{r_i(x)}}{2}\right) \cos \left(\cfrac{\theta_i(x)}{2}\right)\psi(x/\rho_i),\qquad \mbox{for } x\in \mathbf{S}_\delta \mbox{ and }i=1,2,3,
\end{eqnarray*}
where $\psi \in C^\infty_c(\R_+,\R)$ is a smooth cut-off function such that $\psi \equiv 1$ on $[0,1/2]$ and $\psi \equiv 0$ on $[1,+\infty)$.
The positive numbers $\rho_1,\rho_2,\rho_3$ are chosen such that for $\{x \in M_j\;:\; j\neq i\}$ we have $r_i(x) >\rho_i$ . In particular the functions $s_i$ have disjoint supports.
\end{defi}

\begin{figure}[h]
\psfrag{d}[][]{ $\! \!\!\! M_{1}$}
\psfrag{e}[][]{$\,M_2$}
\psfrag{f}[][]{$\! M_{3}$}
\psfrag{g}[][]{ $ 2\pi /3$}
\psfrag{h}[][]{$\!\!2\pi /3$}
\psfrag{i}[][]{$\xi_1$}
\psfrag{j}[][]{$\xi_2$}
\epsfxsize 10cm
\hspace{0cm}
\rotatebox{0}{
\includegraphics[scale=0.4]{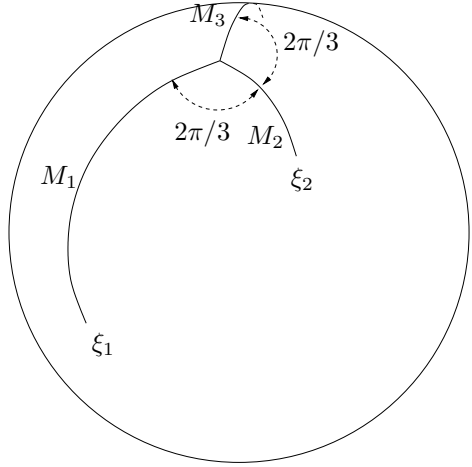}}
\caption{$\mathbf{S}_\delta$}
\label{figure3}
\end{figure}

For $i=1,2,3$, the function $s_i$ defined above belongs to $H^1(\mathbf{S}_\delta)$. Moreover, this function satisfies homogeneous Neumann boundary conditions on $\partial \mathbf{S}_\delta$, we have $\int_{\mathbf{S}_\delta} s_i =0$ and $\Delta s_i$ belongs to $L^2(\mathbf{S}_\delta)$. In fact $\Delta s_i\equiv 0$ on $\{x\;:\; r_i(x)<\rho_i/2 \}$. If $\mathbf{S}_\delta$ were a domain with a smooth boundary, then the quoted properties would imply: $s_i \in H^2(\mathbf{S}_\delta)$. In fact, we have $s_i \in H^s(\mathbf{S}_\delta)$ if and only if $s<3/2$. 

\begin{thm}
\label{thm4}
Let $f \in L^2(\mathbf{S}_\delta)$ such that $\int_{\mathbf{S}_\delta} f =0$, and let $\Sigma \in H^1(\mathbf{S}_\delta)$ solves
\begin{eqnarray*}[\{]
\Delta \Sigma &=&f, \qquad \mbox{in }\: \mathbf{S}_\delta,\\
\partial_n \Sigma &=&0,\qquad \mbox{on } \partial \mathbf{S}_\delta.
\end{eqnarray*} 
Then there exists $\tilde{\Sigma} \in H^2(\mathbf{S}_\delta)$ and $\alpha_1,\alpha_2,\alpha_3 \in \R$ such that
\begin{eqnarray*}
\Sigma &=&\tilde{\Sigma}+\sum_{i=1}^3 \alpha_i s_i.
\end{eqnarray*}
\end{thm}

\section{Necessary conditions}
\label{s2}
Let $(u^\star,K^\star)$, $\delta$ and $\Sigma^\star$ be as in Hypothesis~\ref{hypothese}. 

\noindent
{\bf{$1^{\mbox{st}}$ condition}}

\noindent
Since $(u^\star,K^\star)$ is a global minimizer, the function $u^\star$ belongs to $H^1(B(0,r)\setminus K^\star)$ for every $r>0$, thus  $\Sigma^\star \in H^1(\mathbf{S}_\delta)$ and from Remark~\ref{rmk1.1}, we have:
\begin{eqnarray}
\label{H1}
\Sigma^\star \in H^1(\mathbf{S}_\delta)\setminus \{0\}.
\end{eqnarray} 

\newpage
\noindent
{\bf{$2^{\mbox{nd}}$ condition}}

\noindent
Moreover, from~\eqref{inegalite} with $L=K^\star$, we have for every $r>0$ and every $v$ in $H^1(B(0,r)\setminus K^\star)$ such that $v_{|\partial B(0,r)}=u_{|\partial B(0,r)}$:
\begin{eqnarray*}
\int_{B(0,r)} |\nabla u^\star|^2 &\leq & \int_{B(0,r)} |\nabla v|^2.
\end{eqnarray*}
We deduce that $u^\star$ is harmonic in $\R^3 \setminus K^\star$ and satisfies homogeneous Neumann boundary conditions on $K$. In term of $\Sigma^\star$, the last assertion reads
\begin{eqnarray}
\label{euler}
\left\{ \begin{array}{rcl}
-\Delta \Sigma^\star =&3/4 \Sigma^\star \qquad &\mbox{in~~}~~ \mathbf{S}_\delta,\\
\partial_n \Sigma^\star=&0 \qquad \quad &\mbox{on~~}~~ \partial \mathbf{S}_\delta.
\end{array} \right.
\end{eqnarray}
In particular 
\begin{eqnarray}
\label{moynulle}
\int_{\mathbf{S}_\delta} \Sigma^\star& = & 0. 
\end{eqnarray}

\noindent
{\bf{$3^{\mbox{rd}}$ condition}}

\noindent
By the uniqueness assumption in Hypothesis~\ref{hypothese}, $K^\star$ is unique up to rotation and translation. Thus, at least two of the lengths $\delta_i$ are equal. In the sequel, we will assume without loss of generality that $\mathbf{S}_\delta$ is symmetric with respect to $P$, i.e:
\begin{eqnarray}
\label{sym1}
\delta_2=\delta_3.
\end{eqnarray}
By uniqueness we also have
\begin{eqnarray}
\label{sym2}
&\Sigma^\star(x,y,z)=&\Sigma^\star(-x,y,z)  ,\qquad \forall\, (x,y,z) \in \mathbf{S}_\delta,\\
\label{sym3}
\mbox{or }&  \Sigma^\star(x,y,z)=&-\Sigma^\star(-x,y,z)  ,\qquad \forall\, (x,y,z) \in \mathbf{S}_\delta.
\end{eqnarray}

\noindent
{\bf{$4^{\mbox{th}}$ condition}}

From \eqref{H1},\eqref{euler},\eqref{moynulle}, we may apply Theorem~\ref{thm4} with $f=-3/4 \Sigma^\star$ and $\Sigma=\Sigma^\star$. There exist $\tilde{\Sigma}^\star \in H^2(\mathbf{S}_\delta)$ and $\alpha_1^\star,\alpha_2^\star,\alpha_3^\star \in \R$ such that 
\begin{eqnarray}
\label{decomp}
\Sigma^\star &=&\tilde{\Sigma}^\star + \sum_{i=1}^3 \alpha_i^\star s_i.
\end{eqnarray}
Now let us return in $\R^3$, let $1 \leq i \leq 3$ and let us define a blow-up sequence $(u_k,K_k)_k$ of $(u^\star,K^\star)$ at $\xi_i$.  Let $D_i$ be the line $\R \xi_i$ and $P_i$ be the half plane containing $M_i$ and whose edge is $D_i$. It is clear that the sequence $(K_k)$ converges to $P_i$ locally for the Hausdorff distance. Let us study the $L^1_{loc}$ convergence of the sequence $(u_k)_k$. We denote by $\Pi$ the map $\R^3\setminus\{0\}\rightarrow \mathbf{S}^2$ defined by $\Pi(x):=x/|x|$. Using the above decomposition, we have
\begin{eqnarray*}
u_k(y)&=& t_k^{-1/2} u^\star \left( \xi_i+t_k y\right) 
=  \left| \xi_i/t_k+ y\right|^{1/2} \Sigma^\star\left(\Pi(\xi_i+t_k y)\right).
\end{eqnarray*}
Let $B$ an open ball ball of $\R^3$. Using the decomposition~\eqref{decomp} and the fact that for $j\neq i$, $s_j\equiv 0$ in the neighborhood of $\xi_j$, we have for $t_k$ small enough:
\begin{eqnarray*}
u_k(y)&=& \left| \xi_i/t_k+ y\right|^{1/2} \left( \alpha_i^\star s_i \left(\Pi(\xi_i+t_k y) \right)+ \tilde{\Sigma}^\star\left( \Pi(\xi_i+t_k y)\right)\right),\\
&=& \I_k(y) + \II_k(y).
\end{eqnarray*}
Now we introduce new polar coordinates in $\R^3$: $(R_i(y),\Theta_i(y), z_i(y))$ such that $D_i=\{y\;:\;R_i(y)=0\}$, $P_i=\{y\;:\;\Theta_i(y)=0\}$ and the azimuth is uniquely defined by $z_i(\xi_i)=1$ and $z_i(0)=0$. We have  
\begin{eqnarray*}
\theta_i(\Pi(\xi_i+t_k y ))&=&\Theta_i(y) + \mathcal{O}(t_k),\\
r_i(\Pi(\xi_i+t_k y))&=&t_k R_i(y)+\mathcal{O}(t_k^2),\\
 \left| \xi_i/t_k+ y\right|^{1/2}&=&t_k^{-1/2} + \mathcal{O}(1),
\end{eqnarray*}
uniformely in $y \in B$. Thus, from the definition of $s_i$, we obtain 
\begin{eqnarray*}
\I_k(y)=\alpha_i^\star \sqrt{R_i(y)} \cos \frac{\Theta_i(y)}{2}+\mathcal{O}(\sqrt{t_k}),
\end{eqnarray*}
For the second term, we use the fact that $H^2(\mathbf{S}_\delta)$ is embedded in the H\"older space $C^{0,\gamma}$ for $0<\gamma<1$, in particular choosing $\gamma>1/2$, we easily obtain
\begin{eqnarray*}
\II_k(y)-t_k^{-1/2} \tilde{\Sigma}^\star(\xi_i)&=& o(1),
\end{eqnarray*}
 uniformely in $B$ and we deduce that $u_k-t_k^{-1/2} \tilde{\Sigma}^\star(\xi_i)$ converges to 
\begin{eqnarray*}
u^\star_i(y)&:=& \alpha_i^\star \sqrt{R_i(y)} \cos \frac{\Theta_i(y)}{2},\qquad \mbox{in } L^1_{loc}(\R^3).
\end{eqnarray*}
We now use the following result~\cite{David2}
\begin{thm}
\label{thmglobglob}
Every blow-up limit of a global minimizer is a global minimizer.
\end{thm}
The pair $(u^\star_i,P_i)$ is thus a global minimizer and since we have supposed that the list of global minimizers was closed, the only possibility is that  $(u^\star_i,P_i)$ is a global minimizer of type (v). Consequently, we have
\begin{eqnarray*}
|\alpha_i^\star|&=&\sqrt{2/\pi},\qquad \forall 1\leq i \leq 3.
\end{eqnarray*}
The first consequence of this equality is to exclude the case $\Sigma^\star$ symmetric (Eq. \ref{sym2}). Indeed, in this case, we would have $\alpha_1^\star=0$. Thus $\Sigma^\star$ is antisymmetric, this symmetry implies $\alpha_2^\star=\alpha_3^\star$ and the additional information given by the last equality may be written: 
\begin{eqnarray}
\label{condcoeffs}
|\alpha^\star_1|&=&|\alpha_2^\star|.
\end{eqnarray}  

\noindent
{\bf{$5^{\mbox{th}}$ condition}}
Let $L^2_{0,A}(\Sd)$ and $V^1_A(\Sd)$ be the subspaces of antisymmetric functions in $L^2_0(\Sd)$ and $V^1(\Sd)$. If \eqref{H1} and \eqref{sym3} are true, then, in particular $4/3$ is an eigenvalue of the operator $\Delta_{N,\delta}^{-1}$. Clearly,  $L^2_{0,A}(\Sd)$ is stable by the operator $-\Delta_{N,\delta}^{-1}$. Let us note $-\Delta_{N,\delta,A}^{-1}$ this restriction and $\mu_{a,A}(\delta) \geq \mu_{2,A}(\delta)\geq \cdots >0$ its eigenvalues counting multiplicities. We set $\lambda_{k,A}(\delta):=1/\mu_{k,A}(\delta)$. For $k\geq 1$, we have
\begin{eqnarray}
\label{local}
\lambda_{k,A}(\delta)&:=& \min_{V_k} \max \left\{\int_\Sd |\nabla S^2| \;:\; S\in V_k,\, \int_\Sd S^2 =1\right\},
\end{eqnarray}
where the minimum is taken over all $k$-dimensional subspaces $V_k$ of $V^1_{A}(\Sd)$.

The next results states that $\lambda_{2,A}(\delta)\geq 2$. Consequently,~\eqref{H1},~\eqref{euler},~\eqref{moynulle} and \eqref{sym3} imply that 
\begin{eqnarray}
\label{ppvp}
3/4=\lambda_{1,A}(\delta).
\end{eqnarray}

\begin{prop}
For $\delta \in (0,\pi)^3$, $\lambda_{2,A}(\delta)\geq 2$. 
\end{prop}
\begin{proof}
Let $\delta$ in $(0,\pi)^3$. We have $V^1_A(\mathbf{S}_{(\pi,\pi,\pi)}) \supset V^1(\Sd)$ (where $\mathbf{S}_{(\pi,\pi,\pi)}:=\mathbf{S}^2 \setminus \cup_{i=1}^3 C_i$). Thus, from~\eqref{local},
\begin{eqnarray*}
\lambda_{2,A}(\delta)\geq \lambda_{2,A}(\pi,\pi,\pi),
\end{eqnarray*}
where for $k\geq1$,
\begin{eqnarray*}
\lambda_{k,A}(\pi,\pi,\pi)&:=&\min_{ \stackrel{ V_k\subset V^1(\mathbf{S_{(\pi,\pi,\pi)}}),}{\,\dim V_k =k}} \max \left\{\int_\Sd |\nabla S^2| \;:\; S\in V_k,\, \int_\Sd S^2 =1\right\}.
\end{eqnarray*} 
We have $\lambda_{1,A}(\pi,\pi,\pi)=0$ with associated eigenspace $\R \Sigma_1$ where $\Sigma_1\equiv i$ on $\Omega_i$ for $-1\leq i \leq 1$. It is not difficult to see that there is no other eigenvector in $V^1_A(\mathbf{S}_{(\pi,\pi,\pi)})$ which is locally constant.

\noindent
 Now, let $\Sigma_2 \neq 0$ be an eigenvector associated to $\lambda_{2,A}(\pi,\pi,\pi)$. We split $\Sigma_2$ in $\Sigma_2=S_{-1}+S_0+S_1$, where $\supp S_i \subset \Omega_i$. Let us fix $i$ such that $S_i\not \equiv 0$. This function is a non constant eigenvector of $-\Delta$ restricted to $\Omega_i$ satisfying Neumann boundary conditions, in particular $\int_{\Omega_i} S_i =0$. We set $S:=S_i\circ R^i$, so that $\supp S \subset \bar{\Omega_0}$. We also define $\bar S$ by 
\begin{eqnarray*}
\bar S (x,y,z)&=&S(-x,y,z),\qquad \forall \, (x,y,z) \in \mathbf{S}_{(\pi,\pi,\pi)}.
\end{eqnarray*}  
We have to study two cases.

\noindent
case 1:~~~$S\equiv\bar{S}$. Since $S$ is symmetric, we can define $\Sigma$ in $V^1(\mathbf{S}^2)$ by
 \begin{eqnarray*}
\Sigma(x,y,z) &:=&\left\{
\begin{array}{ll}
S(x,y,z)  &\mbox{ if } (x,y,z) \in  \Omega_0, \\
S(R(x,y,z))&\mbox{ if } (x,y,z) \in  \Omega_{-1},\\
S(R^{-1}(x,y,z))& \mbox{ if } (x,y,z) \in  \Omega_{1}.
\end{array}
\right.
\end{eqnarray*}

\noindent 
case 2:~~~We set $S':=S-\bar{S} \not \equiv 0$. This function is antisymmetric, in particular $S'\equiv 0$ on $P$. In this case, we set:
\begin{eqnarray*}
\Sigma(x,y,z) &:=&\left\{
\begin{array}{ll}
S'(x,y,z)  &\mbox{ if } (x,y,z) \in  \Omega_0, \\
S'(-x',y',z')\quad \mbox{ where } (x',y',z')=R(x,y,z),&\mbox{ if } (x,y,z) \in  \Omega_{-1}\cap \{y\leq 0\},\\
S'(-x',y',z')\quad \mbox{ where } (x',y',z')=R^{-1}(x,y,z),& \mbox{ if } (x,y,z) \in  \Omega_{1}\cap \{y\leq 0\},\\
0& \mbox{otherwise}.
\end{array}
\right.
\end{eqnarray*}

In both cases $\Sigma \in V^1(\mathbf{S}^2)\setminus \{0\}$, satisfies $\int_{\mathbf{S}^2} \Sigma =0$ and 
\begin{eqnarray*}
\int_{\mathbf{S}^2} |\nabla \Sigma|^2 &=&\lambda_{2,A}(\pi,\pi,\pi) \int_{\mathbf{S}^2} \Sigma^2.
\end{eqnarray*}
Thus, from~\eqref{starmorgan}, we conclude that $\lambda_{2,A}(\pi,\pi,\pi) \geq 2$.
\end{proof}

\noindent
{\bf Complete problem. }
We now collect the necessary conditions obtained in this section. If Hypothesis~\ref{hypothese} is true, from~(\ref{H1},\ref{euler},\ref{moynulle},\ref{sym1},\ref{sym3},\ref{condcoeffs} and \ref{ppvp}), then there exists $\delta=(\delta_1,\delta_2,\delta_2) \in (0,\pi)^3$ such that:
\begin{eqnarray}
\label{p1}
3/4=\lambda_{1,A}(\delta)&=& \min \left \{ \int_{\Sd} |\nabla S|^2 \;:\; S\in V^1(\Sd), \,\int_{\Sd} S^2 =1,\, S \mbox{ antisymmetric}\right\}.  
\end{eqnarray}
Moreover, letting
\begin{eqnarray*}
\Sigma(\delta)&\in& \argmin \left \{ \int_{\Sd} |\nabla S|^2 \;:\; S\in V^1(\Sd), \, \int_{\Sd} S^2 =1,\, S \mbox{ antisymmetric}\right\},  
\end{eqnarray*}
then the singular coefficients $\alpha_{1}(\delta),\alpha_{2}(\delta),\alpha_3(\delta)$ such that $\Sigma(\delta)-\sum_{1\leq i \leq 3}\alpha_{i}(\delta) s_i \in H^2(\Sd)$ satisfy
\begin{eqnarray}
\label{p2}
|\alpha_{1}(\delta)|&=&|\alpha_{2}(\delta)|.
\end{eqnarray} 
In the sequel, we give numerical evidences showing that there is no pair $(\delta$, $\Sigma(\delta))$ satisfying both~\eqref{p1} and~\eqref{p2}. The conclusion is that~Hypothesis~\ref{hypothese} is false.

\section{Numerical methods}
\label{s3}
The general method is the following. Let $h>0$, we choose a subdivision  $0=\delta_0^h<\delta_1^h<\cdots<\delta_N^h<\delta_{N+1}^h=\pi$, satisfying $\delta_{k+1}^h-\delta_k^h<h$ for $0\leq k \leq N$. Then, for every $\delta^h:=(\delta_{k_1}^h,\delta_{k_2}^h,\delta_{k_2}^h)$ ($1\leq k_1,k_2 \leq N$), we compute numerical approximations of $\lambda_{1,A}(\delta^h)$ and of the coefficients $\alpha_{i}(\delta^h)$. Finally, we use these values to test the validity of equalities~\eqref{p1} and \eqref{p2}. 

We use a Galerkin method to approximate $\lambda_{1,A}(\delta^h)$. More precisely, we set  
\begin{eqnarray}
\label{p3}
\lambda_{1,A}^h(\delta^h)&:=&\min   \left\{\int_\Sdh |\nabla S^2| \;:\; S\in V^h(\delta^h),\, \int_\Sdh S^2 =1\right\},\\
\label{p4}
\Sigma^h(\delta_h)&\in & \argmin \left\{\int_\Sdh |\nabla S^2| \;:\; S\in V^h(\delta^h),\, \int_\Sdh S^2 =1\right\},
\end{eqnarray} 
where $V^h(\delta^h)$ is a finite dimensional subspace of $V^1(\Sdh)$. This space is chosen great enough such that we may hope that $\lambda_{1,A}^h(\delta^h)$ and $\Sigma^h(\delta_h)$ are close to $\lambda_{1,A}(\delta^h)$ and $\Sigma(\delta_h)$. Typically, $V^h(\delta^h)$ is the space of $P^1$ finite elements constructed on a triangular mesh of $\Sdh$ of size $h$.

\begin{rmk}
We use the same letter ($h$) to denote the step size of the subdivison $\delta_0^h<\cdots<\delta_{N+1}^h$ and the mesh size of the triangular mesh of $\Sdh$. These sizes could be different but they are actually equal in the numerical computations below. 
\end{rmk}

Let $(\Th)_{h>0}$ be a family of regular meshes of $\mathcal{S}^2$ composed of geodesic triangles and with mesh size $h$. We assume that the edges of $\Th$ do not cross the geodesic segments $\mathcal{C}_1$, $\mathcal{C}_2$, $\mathcal{C}_3$. We also assume that $R(\Th)=\Th$ and that $\Th$ is symmetric with respect to $P$. This last symmetry is imposed in order to work with antisymmetric functions.  We choose the subdivision  $0=\delta_0^h<\delta_1^h<\cdots<\delta_N^h<\delta_{N+1}^h=\pi$ such that $(0,\sin \delta_i^h,\cos \delta_i^h)_{0\leq i \leq N+1}$ are the coordinates of the nodes of $\Th$ belonging to $\mathcal{C}_1$.

\noindent
 Let $h>0$. From now on, $\delta^h=(\delta_{k_1}^h,\delta_{k_2}^h,\delta_{k_2}^h)$ and to lighten notations, references to $\delta^h$ will be omitted. Let $(p_i)_{1\leq i \leq M^h}$ be the set of nodes of $\Th$ and let $\mathcal{P}_h$ be the polyhedral domain of vertices $(p_i)_i$ (the boundary of the convex hull generated by $(p_i)_i$). Recall that $\Pi$ is the projection of $\R^3\setminus \{0\}$ on $\mathbf{S}^2$. This map defines a bijection from $\mathcal{P}_h$ onto $\mathcal{S}_2$, let us note $\Pi^{-1}$ its inverse. 

\noindent
Now let $(\bar{\varphi}_i^h)_{1\leq i \leq N^h}$ be the set of continuous functions defined on $\mathcal{P}_h\setminus \cup_{k=1}^3 \Pi^{-1}(M_k)$ such that the restriction of $\bar{\varphi}_i^h$ on each face of $\mathcal{P}_h$ is linear and such that there exists $1\leq j(i) \leq M^h$ such that $\bar{\varphi}^h_i(p_j)=1$ for $j=j(i)$, $0$ otherwise.

\noindent
Finally, for $1\leq i \leq N^h$, we set $\varphi_i^h:=\bar{\varphi}_i^h\circ \Pi^{-1}$ and we define the space of $P^1$ finite elements on $\Sdh$ to be
\begin{eqnarray*}
\bar{W}^h:= \vect \left\{ \varphi^h_i\;:\; 1 \leq i \leq N^h \right\}.
\end{eqnarray*}
And then
\begin{eqnarray*}
\bar{V}^h:= \left \{ \varphi^h \in \bar{W}^h \;:\; \int_{\Sdh} \varphi^h =0\right\}. 
\end{eqnarray*}  
\begin{rmk}
In general, the elements of $\bar{W}^h$ are not continuous across the geodesic segments $M_1,M_2,M_3$. Let us also stress that we have $N^h>M^h$. Indeed, if $p_j$ belongs to $M_1$ and $p_j\neq \xi_1$ then there exists $i_1\neq i_2$, such that $\varphi_{i_1}(p_j)=\varphi_{i_2}(p_j)=1$ and $\supp \varphi_{i_1}\subset \bar{\Omega_{-1}}$,  $\supp \varphi_{i_2}\subset \bar{\Omega_{1}}$. 
\end{rmk}
\begin{rmk}
\label{projection}
The constant functions belong to $\bar{W}^h$ (indeed, $\sum_{i=1}^{N^h} \varphi_i^h \equiv 1$) and $\bar{V}^h$ is the orthogonal of $1$ for the $L^2$ inner product. 
\end{rmk}


Let $f \in L^2_0(\Sdh)$ and $S:=\Delta^{-1}_{N,\delta^h}$. Since $S$ does not necessarily belong to $H^2(\mathbf{S}(\delta^h))$, the classical convergence rates obtained for the approximation by $P^1$ finite elements for a similar problem on a smooth domain are not valid here. In fact, for quasi uniform meshes, there exists $c>0$, such that  
\begin{eqnarray*}
\min_{S^h\in \bar{V}^h} |S^h - S|_{H^1(\Sdh)} &\geq & c \max_i |\alpha_i| \sqrt{h},
\end{eqnarray*}
where $\alpha_1,\alpha_2,\alpha_3$ are the singular coefficients of $S$. This conclusion holds for $S:=\Sigma(\delta_h)$. To overcome this difficulty, we add the singular functions to the space of finite elements.  Namely, we set: 
\begin{eqnarray}
V^h := \bar{V}^h \oplus \vect \{ s_i\;:\; 1 \leq i \leq 3 \}. 
\end{eqnarray}
This method is called Singular Functions Method (see \cite{BDLN3,BDLN2,BDLN1} for a review on such methods). The usual approximation rates (valid for smooth domains) are recovered.
\begin{eqnarray*}
\min_{S^h\in V^h} |S^h - \Sigma|_{H^1(\Sdh)} &\leq & C h. 
\end{eqnarray*}
A classical result (see~\cite{ErnGuermond}, for example) concerning the approximation of the eigenvectors of an elliptic operator by Galerkin methods leads to 
\begin{eqnarray*}
|\lambda_{1,A}^h(\delta^h)-\lambda_{1,A}(\delta^h)| &\leq & C h^2,\\
|\Sigma^h(\delta^h)-\Sigma(\delta^h)|_{H^1(\Sdh)} &\leq & C h,\\
|\Sigma^h(\delta^h)-\Sigma(\delta^h)|_{L^2(\Sdh)} &\leq & C h^2.
\end{eqnarray*}

For the approximation of the singular coefficients $\alpha_i(\delta^h)$, we use an extraction formula of Moussaoui~\cite{Moussaoui} (see also~\cite{Ciarlet}). We first introduce the dual singular functions
\begin{defi}
\label{defi3.1}
With the notations of Definition~\ref{defi2.1}, we define
\begin{eqnarray*}
S_i(x)&:=& \cfrac{1}{2 \tan {\sqrt{r_i(x)}}/{2}}\, \cos \cfrac{\theta_i(x)}{2}\,\psi(x/\rho_i),\qquad \mbox{for } x\in \mathbf{S}_\delta \mbox{ and }i=1,2,3.
\end{eqnarray*}
Now, for $1\leq i \leq 3$, let $\tilde{p}_i \in V^1(\Sd)$ be the variational solution of
\begin{eqnarray}
\label{ptilde}\left\{\begin{array}{rcl}
-\Delta \tilde{p}_i&=&\Delta S_i,\qquad \mbox{in } \Sd,\\
\partial \tilde{p}_i&=& 0,\qquad \mbox{on } \partial \Sd.
\end{array}
\right.
\end{eqnarray}
Finally, we set 
\begin{eqnarray*}
p_i:=S_i+\tilde{p}_i,\qquad \mbox{for } 1\leq i \leq 3.
\end{eqnarray*}
\end{defi}
\begin{rmk}
We have $\Delta S_i\equiv 0$ on $\{x\;:\; \psi(x/\rho_i)=0\}$, so $\Delta S^i$ is smooth and $\tilde{p}_i$ is well defined. The function $p_i$ does not belong to $H^1(\Sd)$ (we only have $p_i \in H^s(\Sd)$ for $s<1/2$). 
\end{rmk}
\begin{thm}[Moussaoui,~\cite{Moussaoui}]
Let $f \in L^2_0(\Sd)$, $\Sigma:=\Delta^{-1}_{N,\delta} f$ and let $\alpha_1,\alpha_2,\alpha_3$ be the singular coefficients of $\Sigma$. Then 
\begin{eqnarray*}
\alpha_i =\cfrac{1}{\pi} \int_{\Sd} p_i(x) f(x) dx. 
\end{eqnarray*}
\end{thm}
In our case, the singular coefficients $\alpha_i(\delta^h)$ are obtained by the formula above with $f:= \lambda_{1,A}^h(\delta^h) \Sigma(\delta^h)$. In order to get  numerical approximations of these coefficients, we first compute an approximation $\tilde{p}^h_i \in V^h$ of the functions $\tilde{p_i}$. We have  
\begin{eqnarray*}
|\tilde{p}_i^h-\tilde{p}_i|_{L^2}&\leq & C h^2,\qquad 1\leq i \leq 3.
\end{eqnarray*}
Then, we set $p_i^h:=S_i+\tilde{p}_i^h$ and finally:
\begin{eqnarray*}
\alpha_i^h(\delta^h)&:=& \lambda^h(\delta^h) \cfrac{1}{\pi} \int_{\Sdh} p_i^h \Sigma_i^h(\delta^h).
\end{eqnarray*}
The numerical convergence rate is given by
\begin{eqnarray*}
|\alpha_i^h(\delta^h)-\alpha_i(\delta^h)|& \leq & C h^2,\qquad 1\leq i \leq 3.
\end{eqnarray*}

Figure~4 represents the error $e(h)$ on the computation of $\lambda_{1,A}(\delta^h)$, $\alpha_1(\delta^h)$ and $\alpha_2(\delta^h)$ for $1/193 \leq h \leq 1/5$ and $\delta^h=(\pi/2,\pi/2,\pi/2)$. The ``exact'' solution is obtained with $h=1/320$.

For the choice of $\psi$ (Definitions~\ref{defi2.1} and~\ref{defi3.1}), it is sufficient to have a $C^2$ function (we have used a piecewise polynomial function). The main obstacle for the accuracy of the numerical computations turned out to be the restriction on $\rho_i$ (Definition~\ref{defi2.1}). If $\delta^h=(\delta_1^h,\delta_2^h,\delta_2^h)$ is such that one of the $\delta_i^h$ is close to $0$ or $\pi$, then we have to choose a very small $\rho_i$. Consequently, the function  $\Delta S_i$ has great values and we need a fine mesh to get an accurate approximation of $\tilde{p}_i$. For this reason, we have worked with this method for $\delta^h \in (0.1,3.04)^3$. For other values of $\delta^h$, we use a method based only on finite elements (without singular functions) described below.

\begin{figure}[h]
\label{figconv}
\epsfxsize 20cm
\hspace{1.5cm}
\psfrag{h}{$h$}
\rotatebox{0}{
\includegraphics[scale=.5]{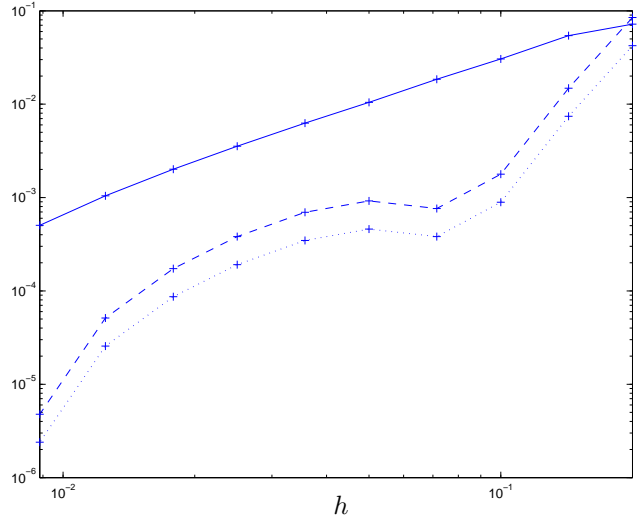}}
\caption{Full line: $|\lambda_{1,A}^h(\delta^h)-\lambda_{1,A}(\delta^h)|$. Dashed line: $|\alpha_{1}^h(\delta^h)-\alpha_{1}(\delta^h)|$. Dotted line: $|\alpha_{2}^h(\delta^h)-\alpha_{2}(\delta^h)|$. $\lambda_{1,A}(\delta^h)\approx 0.795$, $\alpha_1(\delta^h) \approx 0.54$ and $\alpha_2(\delta^h)\approx 0.27$.} 
\end{figure}

Let us fix $\delta \in (0,\pi)^3$. Let $\Sigma^+_i(\delta)$ (resp. $\Sigma^-_i(\delta)$) be the west (resp. east) trace function of $\Sigma(\delta^h)$ on $M_i$. From the definition of $\alpha_i(\delta)$, we have
\begin{eqnarray}
\label{limite}
\lim_{\xi \in M_i \rightarrow \xi_i} \cfrac{\Sigma^+_i(\delta)(\xi)-\Sigma^-_i(\delta)(\xi)}{4 \tan \cfrac{\sqrt{r_i(\xi)}}{2}}&=&\pm \alpha_i(\delta).
\end{eqnarray}
(The sign $\pm$ depends on the orientation choice of Definition~\ref{defi2.1}.)

Now let $h>0$, such that it is possible to set $\delta^h=\delta$. Let $1\leq i \leq 3$, and $\xi^h_0,\cdots,\xi^h_{K_i^h}$ be the nodes of the mesh $\mathcal{T}_h$ belonging to $M_i$. We suppose that the third coordinate of the sequence $(\xi^h)_k$ is decreasing (in particular  $\xi^h_0=(0,0,1)$ and $\xi^h_{K_i}=\xi_i$). We replace $V^h$ by $\bar{V}^h$ in \eqref{p4} to compute an approximation $\Sigma^h_{EF}(\delta^h)$ of $\Sigma(\delta)$ and we define a new approximation of the coefficient $\alpha_i(\delta)$ inspired by~\eqref{limite}. 
\begin{eqnarray}
\label{alphEF}
\alpha_{i,EF}^h(\delta)&:=&\cfrac{\Sigma^{h,+}_{EF,i}(\delta)(\xi_{K_i^h-1})-\Sigma^{h,-}_{EF,i}(\delta)(\xi_{K_i^h-1})}{4 \tan \cfrac{\sqrt{r_i(\xi_{K_i^h-1})}}{2}},
\end{eqnarray}
where $\Sigma^{h,+}_{EF,i}(\delta)$ and $\Sigma^{h,-}_{EF,i}(\delta)$ are the traces of $\Sigma^h_{EF}(\delta^h)$ on $M_i$. In fact, the coefficients $\alpha_{i,EF}^h(\delta)$ do not converge to $\alpha_i(\delta)$, when $h$ goes to 0. However, since we are concerned with the ratio $|\alpha_2^h(\delta^h)|/|\alpha_1^h(\delta^h)|$, it turns out that the method makes sense. During numerical experiments, we have observed that, if we consider a family of quasi-uniform meshes $\mathcal{T}_h$ such that, the family of rescaled meshes $1/h \times(\mathcal{T}_h-\xi_i)$ tends to a fixed mesh $\mathcal{T}_i$ of the plane $\{x\in \R^3\;:\; \xi_i\cdot x =1\}$. Then 
\begin{eqnarray}
\label{train}
\lim_{h \rightarrow 0} \alpha_{i,EF}^h(\delta)&=& c_i \alpha_{i}^h(\delta),
\end{eqnarray}
where $c_i$ is a constant depending on $\mathcal{T}_i$. We did not prove this claim.

In our numerical study, the mesh has the same shape in the neighborhoods of the three points $\xi_1,\xi_2,\xi_3$. Consequently, we have $c_1=c_2=c_3$. Thus we may consider that $|\alpha_{2,EF}^h(\delta)|/|\alpha_{1,EF}^h(\delta)|$ is an approximation of $|\alpha_2(\delta)|/|\alpha_1(\delta)|$.  We have compared the numerical convergence of both methods for this ratio. Since the exact ratio is not known, we have used $|\alpha_{2}^h(\delta)/\alpha_{1}^h(\delta)|$ on a finer mesh ($h \approx 10^{-3}$) to evaluate the error. The numerical convergence results for the two methods are given Figure~5.

\begin{figure}[h]
\epsfxsize 20cm
\psfrag{h}{$h$}
\psfrag{e}{}
\hspace{1cm}
\rotatebox{0}{
\includegraphics[scale=.5]{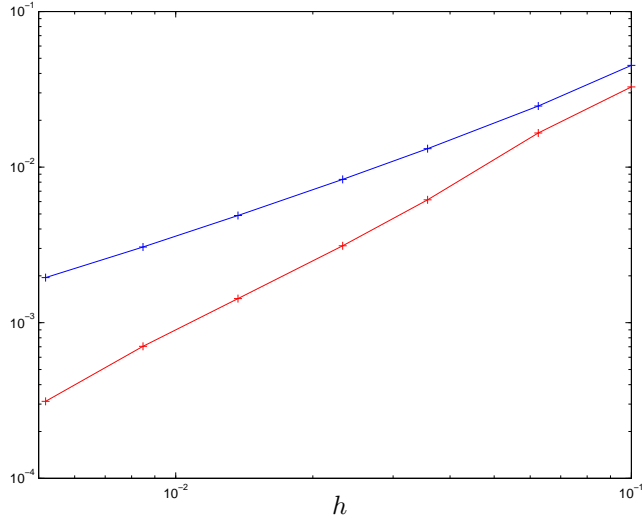}}
\caption{Upper curve: $\left|{|\alpha_{2,\delta,EF}^h|}/{|\alpha_{1,\delta,EF}^h|}-{|\alpha_{2,\delta}|}/{|\alpha_{1,\delta}|}\right|$. Lower curve: $\left|{|\alpha_{2,\delta}^h|}/{|\alpha_{1,\delta}^h|}-{|\alpha_{2,\delta}|}/{|\alpha_{1,\delta}|}\right|$.  For $1/193<h<1/10$, $\delta_1=\pi/6$ and $\delta_2=5\pi/6$.} 
\end{figure}
 The interesting fact is that they do converge to the same limit. We will see in the next section that the zone for which the equalities~\eqref{p1},~\eqref{p2} are the more close to be true is the neighbourhood of ($\delta_1=0$, $\delta_2=\pi$). For this reason, the comparison have been realized for $\delta_1=\pi/6$ and $\delta_2=5\pi/6$.   

\section{Numerical results}
\label{s4}
We have computed $\Sigma^h$ and $\alpha_i^h$ with the numerical method described in the previous section for $h\approx 1/40$. We obtain a curve of approximate solutions of~\eqref{p1}: $\delta_2^h=f^h(\delta^h_1)$ (see Figures~6,7).
\begin{figure}[h]
        \begin{center}
        \begin{minipage}[c]{\linewidth} 
\hspace{.3cm}
        \begin{minipage}[t]{0.50\linewidth}
        \begin{center}
\psfrag{e}[][]{ $\delta_1$}
\psfrag{d}[][]{$\delta_2$}
\psfrag{f}[][]{$\lambda^h_{1,A}(\delta)$}
\centerline{\includegraphics[scale=.5]{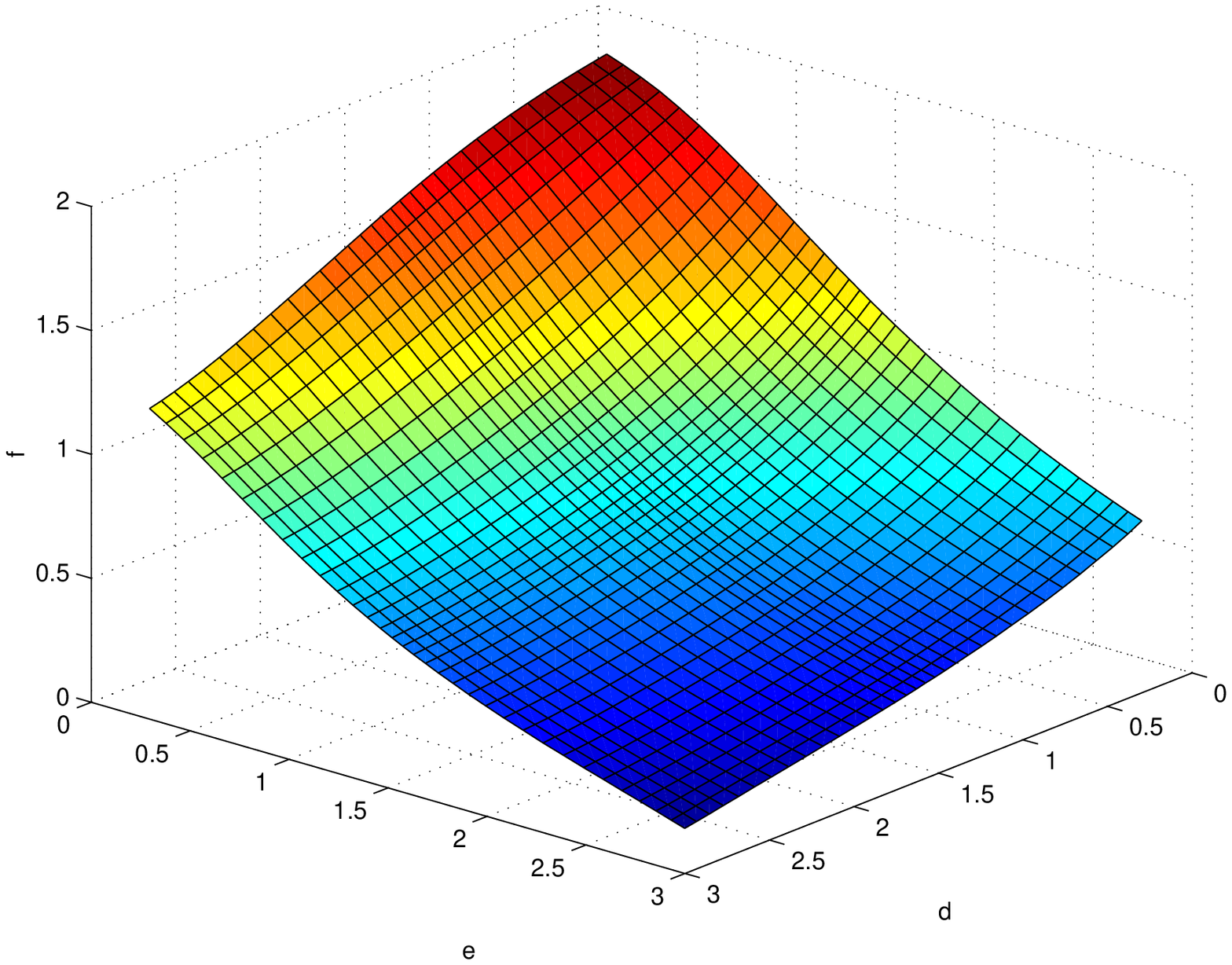}}
\caption{Values of $\lambda_\delta^h$ for $h=1/40$ and $0.14\lesssim \delta_1,\delta_2 \lesssim 3$ .}
 \end{center}
        \end{minipage}
\hspace{-.3cm}
        \begin{minipage}[t]{0.50\linewidth}
        \begin{center}
\psfrag{x}[][]{ $\delta_1=f^h(\delta_2)$}
\psfrag{y}[][]{$\delta_2$}
\vspace{-5.4cm}
\centerline{\includegraphics[scale=.4]{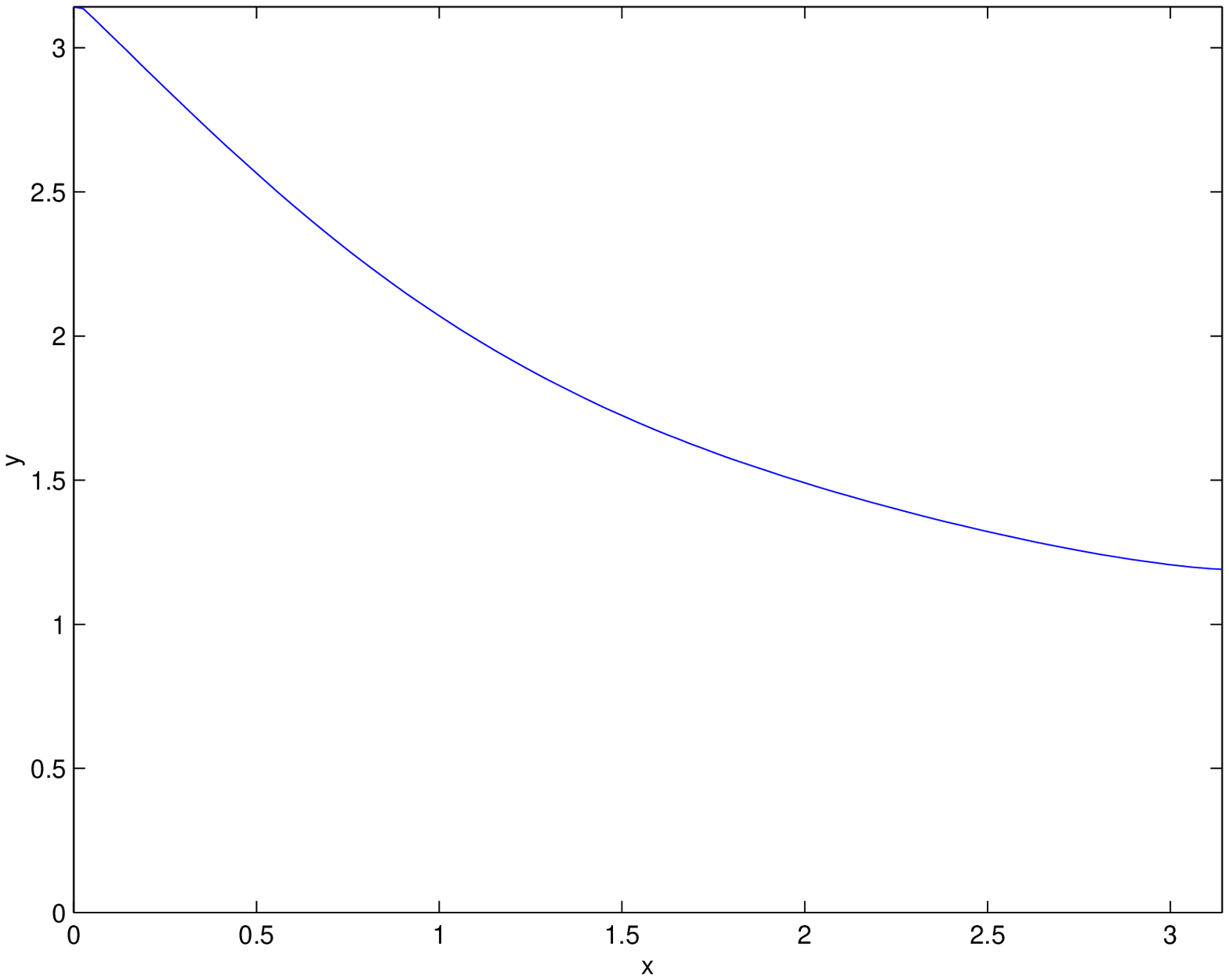}}
\caption{Curve $\delta_1^h:=f^h(\delta_2^h)$, $h=1/160$.}
\end{center}
        \end{minipage}
        \end{minipage}
        \end{center}
\end{figure}
\begin{rmk}
\label{cracktip}
Let us note that for $\delta_1=\pi$ and $\delta_2=0$, the corresponding singular set is a half plane. In this case, a global minimizer do exist: the cracktip (type (v)) and it is natural that $\lim_{(\pi,0)} \lambda_{1,A}^h(\delta_1,\delta_2)=3/4$. Actually, for $(\delta_1,\delta_2)$ close to $(\pi,0)$, we observe that $\Sigma^h(\delta)$ is close to the trace of the cracktip on $\mathbf{S}^2$.

We also have exact values for $(\delta_1,\delta_2)=(0,0)$ and $(\delta_1,\delta_2)=(0,\pi)$ for which $\lambda_{1,A}(0,0)=2$ and $\lambda_{1,A}(0,\pi)=21/16=1.3125$. In the first case, the corresponding eigenvectors are $\R\{(x,y,z)\mapsto y\}$. In the second case, $\S_\delta$ has two connected component, the space of eigenvectors associated to $21/16$ is $\R\Sigma$ where $\Sigma \equiv 0$ on the small connected component and $\Sigma(x,y,z):=(\cos \phi)^{3/4} \sin (3/4 \theta)$ on the big connected component. (The spherical coordinates $(\phi,\theta)$ are defined by $|\theta|<2\pi/3$, $|\phi|<\pi/2$ and $(x,y,z)=(-\cos \phi \cos \theta,-\cos \phi \sin \theta, \sin \phi)$.)
\end{rmk}

In order to check the condition~\eqref{p2}, we compute the approximate coefficients $\alpha_i^h(\delta^h)$ for $h\approx 1/40$ and $0.14\lesssim \delta_1,\delta_2 \lesssim 3$. To complete the study we have computed the alternative approximations of the ratio $|\alpha_{2,\delta}|/|\alpha_{1,\delta}|$ given by $|\alpha^h_{2,EF}(\delta^h)|/|\alpha^h_{1,EF}(\delta^h)|$ for $h\approx 1/160$ on the curve $\{ (f^h(\delta^h_1),\delta_2^h)\;:\; 0<\delta_2^h<\pi\}$ (Figure~8). In both cases, we have 
\begin{eqnarray*}
|\alpha^h_{2}(\delta)|/|\alpha^h_{1}(\delta)| \leq 0.8,
\end{eqnarray*}
for any couple $(\delta_1^h,\delta_2^h)$ of the discretizations. This inequality contradicts~\eqref{p2}.

For both methods, the numerical error is less than $1/100$ (see Figure~5). This numerical error is small compared to the distance between $0.8$ and $1$. We conclude that there is no value $0<\delta_1,\delta_2<\pi$ for which~\eqref{p1} and~\eqref{p2} are both satisfied. Consequently, we are convinced that Hypothesis~\ref{hypothese} was wrong.

\begin{figure}[h]
\psfrag{e}[][]{ $\delta_1$}
\psfrag{d}[][]{$\delta_2$}
\psfrag{f}[][]{$|\alpha_{2}^h(\delta)|/|\alpha_{1}^h(\delta)|$}
\centerline{\includegraphics[scale=.5]{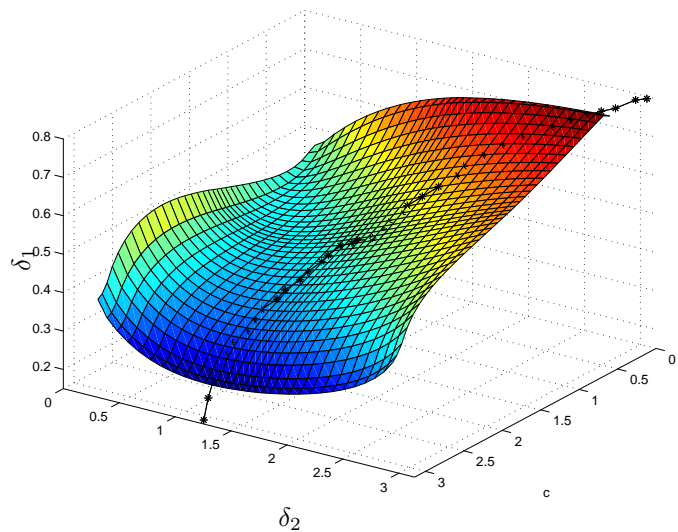}}
\caption{Surface: values of $|\alpha_{2}^h,(\delta^h)|/|\alpha_{1}^h(\delta^h)|$ for $h\approx 1/40$. Black line: values of $|\alpha_{2,EF}^h(\delta^h)|/|\alpha_{1,EF}^h(\delta^h)|$ for $h=1/250$.}
\end{figure}

\begin{rmk}
\label{rm2}
We observe that $\alpha_1(\delta)$ converges to $0$ when $\delta_1$ tends to $0$ so that solutions of~\eqref{p2} do exist, but in this case $\delta_1<f(\delta_2)$ and $\lambda_{1,A}(\delta)>3/4$.
\end{rmk}

\begin{rmk}
Since no solution has been found, we have removed the symmetry condition $\delta_2=\delta_3$. Again we do not find any non zero eigenvector $\Sigma=-4/3 \Delta_{N,\delta}^{-1}\Sigma$ whose singular coefficients satisfy $|\alpha_1|=|\alpha_2|=|\alpha_3|$.
\end{rmk}

\section{Conclusion}
\label{s5}
The above numerical experiments show that Hypothesis~\ref{hypothese} is certainly wrong. The first consequence is that we still don't know the shape of the singular set of a minimizer in the situation of Figures~1,2. One possibility is that the true singular set is topologically equivalent to the one of Figure~2 but with edges tangent to $\gamma$ at $\gamma(1)$ (see Figure~9 below).    

For the moment, this new hypothesis is a conjecture. If it were true, the singular set of a blow-up limit at $\gamma(1)$ would be a half plane and one may expect that the associated global minimizer would be a cracktip (type (v)). In this case there is no need to add a new type of global minimizers to the existing list in order to explain Figure~9. However one may wonder if there exists a global minimizer whose singular set is locally diffeomorphic to the one of Figure~9. Such a global minimizer would not be blow-up invariant. 

Another consequence of this negative result is that taking blow-up limit at $\gamma(1)$, we cannot discrimate a surface with a smooth boundary and the surface above. Thus it seems now more difficult to use the information on global minimizers to deduce some regularity for the singular sets of minimizers.    

\begin{figure}[h]
\label{bourdin}
\psfrag{g}[][]{$\gamma$}
\psfrag{a}[][]{$\!\!\!\gamma(1)$}
\centerline{ \rotatebox{270}{\includegraphics[scale=.4]{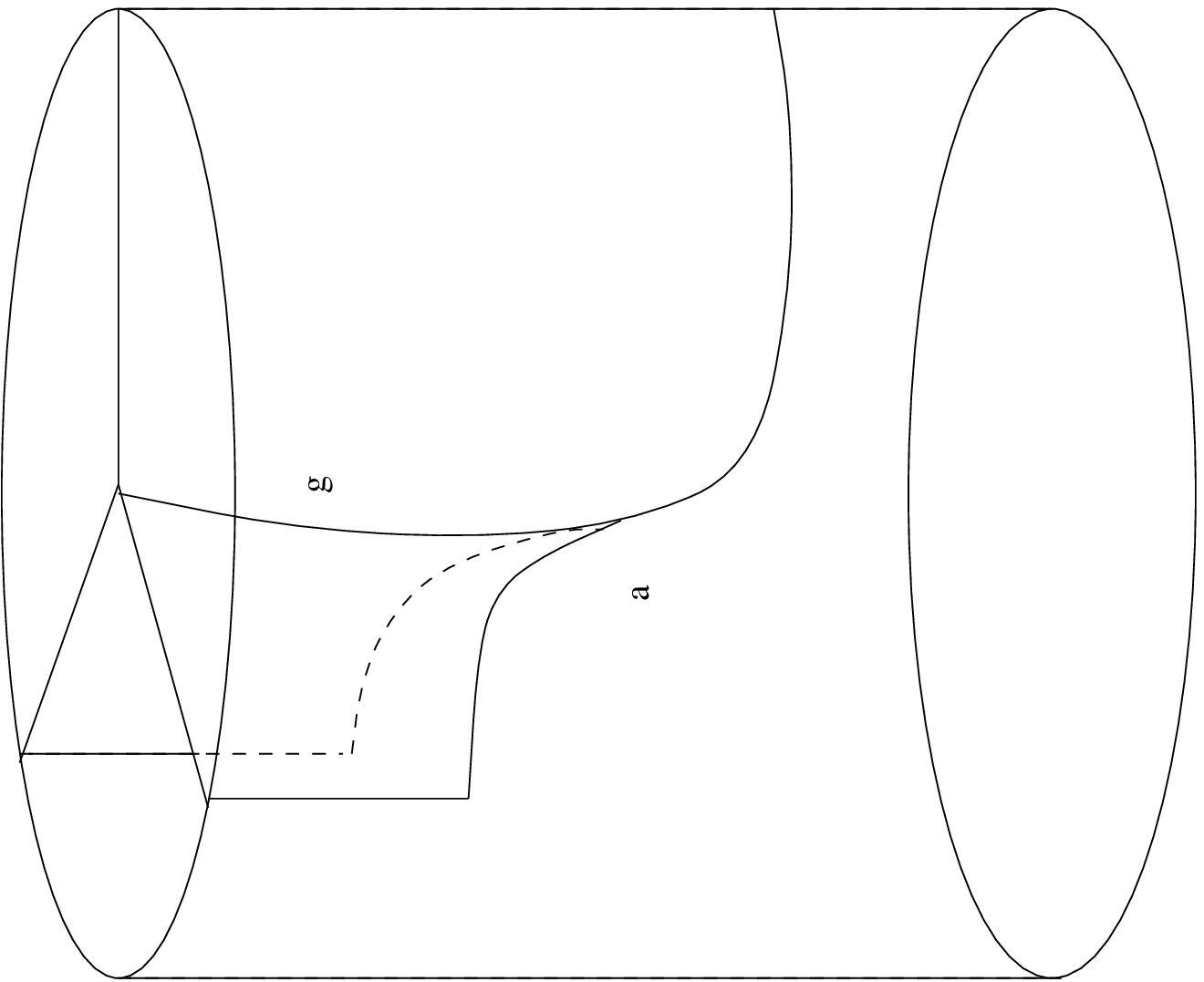}}}
\caption{}
\end{figure}

\bibliographystyle{plain} 
\bibliography{bib} 
 
\end{document}